\documentclass[11pt]{article}
\usepackage[mathscr]{eucal}
\usepackage{amsfonts}
\usepackage{amsmath}
\usepackage{amsthm}
\usepackage{amssymb}
\usepackage{amscd}
\usepackage{latexsym}

\theoremstyle{plain}
\newtheorem{theorem}[equation]{Theorem}

\theoremstyle{definition}
\newtheorem{remark}[equation]{Remark}

\newcommand{\C}{{\mathbb{C}}}
\newcommand{\diag}{{\mathrm{diag}}}
\newcommand{\gl}{{\mathfrak{gl}}}
\newcommand{\h}{{\mathfrak{h}}}
\newcommand{\R}{{\mathcal{R}}}
\renewcommand{\L}{{\mathcal{L}}}
\newcommand{\Z}{{\mathcal{Z}}}
\def\s{{\mathbb{S}}}

\begin{document}

\centerline{\Large {\bf{Chevalley restriction theorem 
for the cyclic quiver}}}

\bigskip

\centerline{\large {\sc Wee Liang Gan}}

\bigskip

\begin{abstract}
We prove a Chevalley restriction theorem and
its double analogue for the cyclic quiver.
\end{abstract}

\bigskip

The aim of this paper is to prove a Chevalley restriction 
theorem and its double analogue for the cyclic quiver.
When the quiver is of type $\widehat{A}_0$, we recover the 
results for $\gl_n$. The proof of our Chevalley restriction
theorem is similar to the proof for $\gl_n$; however, the
proof of the double analogue uses a theorem of 
Crawley-Boevey on decomposition of quiver varieties.
The double analogue is the limiting case of an isomorphism
between a Calogero-Moser space and the center of a symplectic
reflection algebra proved by Etingof and Ginzburg. 
It is also the associated graded version of a conjectural 
Harish-Chandra isomorphism for the cyclic quiver.

We now introduce our notations.
Let $Q$ be the cyclic quiver with $m$ vertices. 
Let $\delta = (1, \ldots, 1)$ be the minimal positive imaginary
root. Let $\R_{n} = \mathrm{Rep}(Q, n\delta)$ be the space of
representations of $Q$ with dimension vector $n\delta$. Thus, 
\[ \R_{n} = \underbrace{\gl_{n} \times \cdots \times \gl_{n}}_{m}.\]
Next, let $\h$ be the subspace of diagonal matrices in $\gl_{n}$, and
let \[ \L_n = \{ (z, \ldots, z)\in \R_{n}\ |\ z \in \h \}.\]
Note that $\L_n$ is a $n$ dimensional subspace of $\R_{n}$.
Let 
\[ G_n = \underbrace{GL_{n} \times \cdots \times GL_{n}}_{m}.\] 
An element $(g_{1}, \ldots, g_{m}) \in G_n$ acts on an element 
$(x_{1}, \ldots, x_{m}) \in \R_{n}$, giving 
\[ (g_{2}^{-1}x_{1}g_{1},\ g_{3}^{-1}x_{2}g_{2},\ \ldots,
\ g_{1}^{-1}x_{m}g_{m}).\]
Let $\s_n$ be the symmetric group on $n$ letters, 
which we will also regard as the subgroup of permutation matrices
in $GL_n$. Finally, let 
$$W_n = \s_{n} \ltimes (\mathbb{Z}/m\mathbb{Z})^{n}.$$
We have $W_n \hookrightarrow G_n$ via
\begin{align*}
(\sigma , \zeta_{1}, & \ldots, \zeta_{n}) \mapsto \\
& ( \sigma\cdot \diag(1, \ldots, 1), \sigma\cdot \diag(\zeta_{1}, 
\ldots, \zeta_{n}), \cdots, \sigma\cdot \diag(\zeta_{1}^{m-1},
\ldots, \zeta_{n}^{m-1})),
\end{align*}
where $\diag(\ldots)$ denotes the diagonal matrix with the 
indicated entries.
Hence, $W_n$ acts on $\R_n$.
Observe that the action of $W_n$ on $\L_n$ is stable.
We remark that $W_n$ is the complex reflection group 
of type $G(m,1,n)$ and $\L_n$ is its reflection representation.

\begin{theorem}[Chevalley restriction] \label{cr}
Restriction of functions from $\R_{n}$ to $\L_n$ gives an 
isomorphism \[ \rho: \mathbb{C}[\R_{n}]^{G_n}
\stackrel{\sim}{\longrightarrow} \mathbb{C}[\L_n]^{W_n}.\]
\end{theorem}

\begin{proof}
{\it Surjectivity of $\rho$}:
Write an element in $\R_{n}$ as
$(x_{1}, \ldots, x_{m})$ and an element in $\L_n$ as
$$(\diag(z_{1}, \ldots, z_{n}), \ldots).$$
Note that $\mathbb{C}[\L_n]^{W_n}$
is a polynomial algebra generated by the elementary symmetric polynomials
in $z_{1}^{m}, \ldots, z_{n}^{m}$. The homomorphism $\rho$ takes
the coefficients of the characteristic polynomial of
$x_{m}x_{m-1}\cdots x_{1}$ to the elementary symmetric polynomials in
$z_{1}^{m}, \ldots, z_{n}^{m}$.
This proved that $\rho$ is surjective.

{\it Injectivity in the $n=1$ case}: 
Call an element $(x_{1}, \ldots, x_{m}) \in
\R_{1}$ generic if $x_{1}\cdots x_{m} \neq 0$. Observe that
the set of generic elements are Zariski open dense in $\R_{1}$.
Moreover, it is easy to see that in this case,
two generic elements $(x_{1}, \ldots, x_{m})$ and
$(x'_{1}, \ldots, x'_{m})$ are in the same $G_1$-orbit iff
$x_{1}\cdots x_{m} = x'_{1}\cdots x'_{m}$. In particular, 
$\L_1$ intersects every generic orbit.
Hence, $\rho$ is injective.

{\it Injectivity in the general case}: 
Call an element $(x_{1}, \ldots, x_{m})$ in $\R_n$
generic if $x_{m}x_{m-1}\cdots x_{1}$ has pairwise distinct
nonzero eigenvalues. Denote the subset of generic elements in $\R_n$
by $\R'_n$, and let $\L'_n=\L_n\cap\R'_n$.
Observe that $\R'_n$ and $\L'_n$ are, respectively, Zariski open
dense in $\R_n$ and $\L_n$. Moreover, $\R'_n$ is $G_n$-stable and
$\L'_n$ is $W_n$-stable. The injectivity of $\rho$ follows from the 
$n=1$ case and the following claim.

\emph{Claim}: If $(x_{1}, \ldots, x_{m})\in \R'_n$, then
it can be diagonalized, i.e. $G_n$-conjugated to an element in 
$$ \underbrace{R_{1} \times \cdots 
\times R_{1}}_{n} = \underbrace{\h \times \cdots \times \h}_{m}. $$

{\it Proof of Claim}: By our assumption, $x_1, \ldots, x_m$ are
invertible matrices. Moreover, there exists an invertible 
matrix $g$ such that $g^{-1}x_{m}x_{m-1}\cdots x_{1}g$ is diagonal.
Then, using $$(g, x_1g, x_2x_1g, \ldots, x_{m-1}\cdots x_1g)\in G_n,$$
we can conjugate $(x_1, \ldots, x_m)$ to 
$$(1, \ldots, 1, g^{-1}x_{m}x_{m-1}\cdots x_{1}g).$$
This proved the claim, and hence the theorem.
\end{proof}

\begin{remark}
The Jacobian of the morphism $\L_n \rightarrow \L_n/W_n$
at a point
$$(z, \ldots)=(\diag(z_{1}, \ldots, z_{n}), \ldots) \in \L_n$$
is, up to a nonzero constant, equal to
$$(z_1\cdots z_n)^{m-1}\prod_{i<j}(z_i^m-z_j^m).$$
Thus, $\L'_n$ is the set of points where
the Jacobian is nonzero. 
\end{remark}

We now proceed to the double analogue of Theorem \ref{cr}.
Let $\Z_n$ be the zero set of the moment map of the $G_n$-action on
$T^{*}\R_{n} = {\mathrm{Rep}}(\overline{Q}, n\delta)$, 
where $\overline{Q}$ is the double quiver of $Q$. 
Write an element in ${\mathrm{Rep}}(\overline{Q}, n\delta)$ as 
$$(x_{1}, \ldots, x_{m}, y_{1}, \ldots, y_{m})
\in \underbrace{\gl_n\times \cdots \times \gl_n}_{2m}.$$ 
Here, the arrow for $y_{i}$ is opposite to the arrow for $x_{i}$.
In explicit terms, $\Z_n$ is defined by the moment map equations
\[ y_{1}x_{1} - x_{m}y_{m} = 0,\ y_{2}x_{2} - x_{1}y_{1} = 0,\ \ldots.\]
The action of an element $(g_{1}, \ldots, g_{m})\in G_n$ on
${\mathrm{Rep}}(\overline{Q}, n\delta)$ 
is given by the formula
\[ (g_{2}^{-1}x_{1}g_{1},\ \ldots, g_{1}^{-1}x_{m}g_{m}
,\ g_{1}^{-1}y_{1}g_{2},\ \ldots,\ g_{m}^{-1}y_{m}g_{1}).\]
Note that $\Z_n$ is stable under the $G_n$-action, and
$\L_n \times \L_n \subset \Z_n$.

\begin{theorem}[Double analogue] \label{da}
Restriction of functions from $\Z_n$ to $\L_n \times \L_n$ gives an
isomorphism 
\[ \phi: \mathbb{C}[\Z_n]^{G_n} \stackrel{\sim}{\longrightarrow}
\mathbb{C}[\L_n \times \L_n]^{W_n}.\]
\end{theorem}

\begin{proof}
{\it Surjectivity of $\phi$}: Write an element of $\L_n \times \L_n$ as
\[ (\diag(z_{1}, \ldots, z_{n}), \ldots, \diag(z'_{1}, \ldots, z'_{n}),
\ldots). \] By a result of Weyl \cite{We},
the algebra $\mathbb{C}[\L_n \times \L_n]^{W_n}$ is generated by
$$z_{1}^{r}{z'_{1}}^{s} + \cdots + z_{n}^{r}{z'_{n}}^{s},$$ where
$r, s \geq 0$ and $r-s$ is divisible by $m$. The homomorphism $\phi$ takes
\[ {\mathrm{Tr}}
( \underbrace{y_{1}\cdots y_{m}y_{1}\cdots y_{m}\cdots y_{1}
\cdots y_{j}}_{s}
\underbrace{x_{j}\cdots x_{1}\cdots x_{m}\cdots x_{1}
x_{m}\cdots x_{1}}_{r}) \]
to $z_{1}^{r}{z'_{1}}^{s} + \cdots + z_{n}^{r}{z'_{n}}^{s}$.
Hence, $\phi$ is surjective.

{\it Injectivity in the $n=1$ case}:
Suppose $x_{1}\cdots x_{m} = d^{m} \neq 0$ for some $d$. Let
$g_{i} = x_{1}\cdots x_{i-1}/d^{i-1}$. Then
$g_{i+1}^{-1}x_{i}g_{i} = d$ and
$g_{i}^{-1}y_{i}g_{i+1} = x_{i}y_{i}/d$.
But $x_{1}y_{1} = x_{2}y_{2} = \cdots x_{m}y_{m}$
by the moment map equations. Hence, the $G_1$-saturation of
$\L_1\times \L_1$ is Zariski dense in $\Z_1$.
It follows that $\phi$ is injective.

{\it Injectivity in the general case}:
Observe that
$$\underbrace{\Z_1\times\cdots\times\Z_1}_{n} 
= \big(\underbrace{\h\times\cdots\times\h}_{m}\big)\cap\Z_n
\hookrightarrow \Z_n.$$
This inclusion induces the map $f$ in the following 
commutative diagram.
\[ \begin{CD}
\C[\Z_n]^{G_n} @>f>> \big((\C[\Z_1]^{G_1})^{\otimes n}\big)^{\s_n} \\
@V{\phi}VV @VV{\psi}V \\
\C[\L_n\times \L_n]^{W_n} @>>> \big((\C[\L_1\times 
\L_1]^{W_1})^{\otimes n}\big)^{\s_n}
\end{CD} \]
In this diagram, the map $f$ is injective by \cite[Theorem 3.4]{CB},
and the map $\psi$ is injective by the $n=1$ case which we have proved.
Hence, $\phi$ must be injective.
\end{proof}

When $Q$ is of type $\widehat{A}_0$, i.e. the $\gl_n$-case,
the injectivity of $\phi$ was due to Gerstenhaber \cite{Ge},
cf. \cite{Ri}. 

\begin{remark}
Let $\Gamma$ be any finite subgroup of $SL_2(\C)$ and let
${\mathbf{\Gamma}}_n = \s_n \ltimes \Gamma^n$.
In \cite{EG}, Etingof and Ginzburg defined the Calogero-Moser space
$\mathcal{M}_{\Gamma,n,c}$ associated to $\Gamma$ and a 
parameter $c$. They proved that, for generic $c$, 
there is an isomorphism $\C[\mathcal{M}_{\Gamma,n,c}]
\stackrel{\sim}{\rightarrow} \mathsf{Z}_{0,c}({\mathbf{\Gamma}}_n)$,
where $\mathsf{Z}_{0,c}({\mathbf{\Gamma}}_n)$ is the center
of the symplectic reflection algebra associated to the group
${\mathbf{\Gamma}}_n$, see \cite[Theorem 11.16]{EG}.
The isomorphism $\phi$ in Theorem \ref{da} is 
the $c=0$ case of their isomorphism when $\Gamma=
\mathbb{Z}/m\mathbb{Z}$;
it is also the associated graded version of the isomorphism
in \cite[Conjecture 11.22]{EG} for this $\Gamma$. 
\end{remark}

{\it Acknowledgment}. I thank Victor Ginzburg for suggesting 
the problem and for helpful discussions.

\footnotesize{
}

Department of Mathematics, Massachusetts Institute of Technology,
Cambridge, MA 02139, U.S.A.;\\
\hphantom{x}\quad\, {\it E-mail address}: {\tt wlgan@math.mit.edu}

\end{document}